\documentclass{article}

\usepackage{citesort}
\usepackage{graphicx}
\usepackage{amssymb}
\usepackage{amsmath}

\newcommand{\D}{\mathrm{d}}  

\def\XXint#1#2#3{{\setbox0=\hbox{$#1{#2#3}{\int}$}
     \vcenter{\hbox{$#2#3$}}\kern-.5\wd0}}

\begin{document}

\title{Numerical quadratures for near-singular and near-hypersingular
  integrals in boundary element methods}

\author{Michael Carley$^*$\\ Department of Mechanical
  Engineering, University of Bath, Bath}


\maketitle

\begin{abstract}
  A method of deriving quadrature rules has been developed which gives
  nodes and weights for a Gaussian-type rule which integrates
  functions of the form:
  \[
  f(x,y,t) = \frac{a(x,y,t)}{(x-t)^{2}+y^{2}} +  
  \frac{b(x,y,t)}{[(x-t)^{2}+y^{2}]^{1/2}} + 
  c(x,y,t)\log[(x-t)^{2}+y^{2}]^{1/2} + d(x,y,t),
  \]
  without having to explicitly analyze the singularities of $f(x,y,t)$
  or separate it into its components. The method extends previous work
  on a similar technique for the evaluation of Cauchy principal value
  or Hadamard finite part integrals, in the case when $y\equiv0$. The
  method is tested by evaluating standard reference integrals and its
  error is found to be comparable to machine precision in the best
  case. 
\end{abstract}

\bibliographystyle{plain}

\section{Introduction}
\label{sec:intro}

An important part of the application of the boundary element method
(BEM) to physical problems is the calculation of the field, from the
solution on the boundary. This requires the evaluation of integrals
which contain `almost-singular' integrands which are not properly
handled by standard Gaussian quadratures. For example, if we consider
a two-dimensional potential problem, such as the Laplace or Helmholtz
equations: 
\begin{subequations}
    \begin{eqnarray}
    \nabla^{2}\phi &= 0,\\
    \nabla^{2}\phi + k^{2}\phi &= 0,
  \end{eqnarray}
\end{subequations}
where $k$ is the wavenumber, the potential $\phi$ at some point in the
field will be:
\begin{equation}
  \label{equ:potential}
  \phi(\mathbf{x}) = \int_{\Gamma} 
  \frac{\partial\phi}{\partial n}G(\mathbf{x},\mathbf{x}_{1})
  -
  \frac{\partial G(\mathbf{x},\mathbf{x}_{1})}{\partial n}\phi\,\D \Gamma,
\end{equation}
where $G(\mathbf{x},\mathbf{x}_{1})$ is the Green's function for the
problem, $\Gamma$ is the boundary of the domain with normal $n$ and
variables of integration are indicated by subscript $1$. In the case
of the Helmholtz equation, the singular behaviour of the Green's
function will be related to the Green's function of the corresponding
Laplace equation~\cite[for example]{dawson95}. The Green's functions
for the Laplace equation will have a logarithmic singularity in the
planar and axisymmetric case (where the Green's function is
proportional to an elliptic integral) and also in the case of an
asymmetric problem in an axisymmetric
domain~\cite{carley06,bjorkberg-kristensson87}.

\begin{figure}
  \centering
  \includegraphics{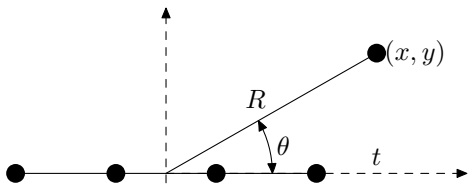}
  \caption{Geometry of boundary element: collocation points indicated
    by filled circles, field point position $x=R\cos\theta$,
    $y=R\sin\theta$.}
  \label{fig:bem:geometry}
\end{figure}

The Green's function for the two-dimensional Laplace equation is
$G=\log|\mathbf{x}-\mathbf{x}_{1}|$. If the boundary integral problem
has been solved using the BEM, the evaluation of the potential in the
field gives rise, on each element, to an integral of the form:
\begin{equation}
  \label{equ:int}
  I^{(0)}(x,y) = \int_{-1}^{1} (\log|(x-t)^{2}+y^{2}|+f(x,y,t))L(t)\,\D t,
\end{equation}
where $t$ is the local coordinate on the element, $L(t)$ is a shape
function, typically a polynomial, and $x$ and $y$ are the field point
coordinates in the local coordinate system,
figure~\ref{fig:bem:geometry}. This integral is non-singular but
suffers from an `offstage singularity'~\cite{acton90} if the field
point is near the element, so that the argument of the logarithm
becomes small.

If it is required to determine the first or second derivatives of the
potential, for example, to determine a velocity or velocity gradient
in a fluid-dynamical problem, integrals of the form:
\begin{subequations}
  \label{equ:int:deriv}
  \begin{eqnarray}
  \label{equ:int:1}
  I^{(1)}(x,y) &= \int_{-1}^{1} \frac{L(t)}{[(x-t)^{2}+y^{2}]^{1/2}} +
  L(t)g(x,y,t)\,\D t,\\
  \label{equ:int:2}
  I^{(2)}(x,y) &= \int_{-1}^{1} \frac{L(t)}{(x-t)^{2}+y^{2}} +
  L(t)h(x,y,t)\,\D t,
  \end{eqnarray}
\end{subequations}
arise. By analogy with the case when $y\equiv0$, we refer to these
integrals as `near-singular' and `near-hypersingular' respectively.
When $y\equiv0$, the field point lies on the boundary and the
integrals must be treated as Cauchy principal
values~\cite[page~37]{lighthill58} or as a Hadamard
finite-part~\cite[page~31]{lighthill58}. In this case, a number of
procedures exist for the accurate evaluation of the
integrals~\cite[for example]{brandao87,monegato87,monegato94}

A recently-developed method~\cite{kolm-rokhlin01,carley07a} develops
quadrature rules which can be used as `plug-in' replacements for
standard Gaussian quadratures in cases where the integrand contains
singularities of the form of equations~(\ref{equ:int})
and~(\ref{equ:int:deriv}). An important feature of these rules is that
they require no analysis of the integrand before use. When the
potential is differentiated, yielding equation~(\ref{equ:int:1}), the
function $g(x,y,t)$ will usually contain logarithmic singularities
which must be properly accounted for, as well as the leading
singularity. In many problems, the Green's function will be
complicated and the explicit identification of the singularities will
be time-consuming. It is preferable to have a quadrature rule which
correctly integrates all of the singularities present without
requiring a detailed analysis. The previously published
method~\cite{kolm-rokhlin01,carley07a}, generates such rules for
singular and hypersingular integrals. This paper extends this method
to `near-singular' and `near-hypersingular' integrals, giving a
simple, easily implemented technique producing quadrature rules which
are direct replacements for Gaussian quadrature in BEM codes.

\section{Quadrature rules for near-singular integrals}
\label{sec:method}

Before developing the quadrature rules needed for near-singular
integrals, it is worth examining the reason for the breakdown of
standard Gaussian quadrature. This will also give an indication of
when specialized rules should not be used and standard quadratures are
better. Consider the integral:
\begin{eqnarray}
  \label{equ:int:sample}
  I &= \int_{-1}^{1} \frac{1}{[(x-t)^{2}+y^{2}]^{1/2}}\,\D t,\\
  &= \frac{1}{R}\int_{-1}^{1} \frac{1}{((t/R)^{2}-2t\cos\theta/R +
    1)^{1/2}} \,\D t,\nonumber
\end{eqnarray}
where $R^{2}=x^{2}+y^{2}$ and $\theta=\tan^{-1}y/x$. This can be
expanded using the generating function for Legendre
polynomials~\cite[8.921]{gradshteyn-ryzhik80}:
\begin{eqnarray}
  \label{equ:generation}
  \frac{1}{(1-2zt+t^{2})^{1/2}} = 
  &\sum_{k=0}^{\infty} t^{k}P_{k}(z), & |t| < \min|z\pm(z^{2}-1)^{1/2}|,\\
  &\sum_{k=0}^{\infty} t^{-(k+1)}P_{k}(z), & |t| >
      \max|z\pm(z^{2}-1)^{1/2}|,
\end{eqnarray}
to yield:
\begin{subequations}
  \label{equ:expansion}
  \begin{eqnarray}
    \label{equ:expansion:1}
    \frac{1}{[(x-t)^{2}+y^{2}]^{1/2}} &=
    \sum_{k=0}^{\infty} \frac{t^{k}}{R^{k+1}}P_{k}(\cos\theta),\quad t<R\\
    \label{equ:expansion:2}
    &=\sum_{k=0}^{\infty}
    \frac{R^{k}}{t^{k+1}}P_{k}(\cos\theta),\quad t>R. 
  \end{eqnarray}
\end{subequations}
An $N$-point Gaussian quadrature rule integrates exactly polynomials
up to order $2N-1$. When $R\gg1$, the integrand is well approximated
by equation~(\ref{equ:expansion:1}) and the estimate of $I$ returned
by Gauss-Legendre quadrature:
\begin{equation}
  \label{equ:gauss:estimate}
  I \approx
  2\sum_{k=0}^{N}\frac{1}{(2k+1)R^{2k+1}}P_{2k}(\cos\theta),
\end{equation}
will be accurate. If, however, $R<1$, part of the integrand will be
given by the inverse power series of equation~(\ref{equ:expansion:2})
which cannot be correctly integrated by the standard Gaussian
quadrature. Likewise, even if $R>1$ but is not large enough to make
the terms of equation~(\ref{equ:expansion:1}) decay fast enough, there
will be a large error in the estimate of $I$. Looking ahead to the
results presented in \S\ref{sec:tests}, it is expected that standard
Gauss-Legendre quadrature will give good results for large $R$ and/or
in cases where $P_{n}(\cos\theta)$ is small. Otherwise, a specialized
rule will be necessary.

\subsection{Evaluation of quadrature rules}
\label{sec:rules}

The algorithm to be developed gives an $N$-point rule which integrates
a function of the form:
\begin{equation}
  \label{equ:quad:func2}
  f(x,y,t) = \frac{a(x,y,t)}{(x-t)^{2}+y^{2}} +  
  \frac{b(x,y,t)}{[(x-t)^{2}+y^{2}]^{1/2}} + 
  c(x,y,t)\log[(x-t)^{2}+y^{2}]^{1/2} + d(x,y,t),
\end{equation}
where $a$, $b$, $c$ and $d$ are taken to be polynomials of order up to
$M$ and the integral
\begin{equation}
  \label{equ:integral}
    I(x,y) = \int_{-1}^{1}f(x,y,t)\,\D t 
    \approx \sum_{i=0}^{N} w_{i}f(x,y,t_{i}),
\end{equation}
where $t_{i}$ are the quadrature points of an $N$-point Gauss-Legendre
quadrature and the weights $w_{i}$ are to be determined. A
previously-developed algorithm~\cite{kolm-rokhlin01,carley07a} for the
computation of quadrature rules gives a method for the evaluation of
$w_{i}$ when $y\equiv0$. The approach is conceptually simple---the
weights are found as the solution to the system of equations:
\begin{equation}
  \sum_{j}[\psi_{ij}]w_{j} = m_{i},\quad i=1,\dots,4M,
\end{equation}
where $\psi_{ij}$ are the weighted Legendre polynomials at the
quadrature points $t_{j}$:
\begin{equation}
  \psi_{ij} =
  \left\{
    \begin{array}{ll}
      P_{i-1}(t_{j}) & 1 \leq i \leq M,\\
      P_{i-M-1}(t_{j})\log [(x-t_{j})^{2}+y^{2}]^{1/2} & M+1 \leq i \leq 2M,\\
      P_{i-2M-1}(t_{j})/[(x-t_{j})^{2}+y^{2}]^{1/2} & 2M+1 \leq i \leq 3M,\\
      P_{i-3M-1}(t_{j})[(x-t_{j})^{2}+y^{2}] & 3M+1 \leq i \leq 4M.
    \end{array}
  \right.
\end{equation}
and the moments $m_{i}$ are the integrals of $\psi_{i}$
\begin{equation}
  \label{equ:moments}
  m_{i} = 
  \left\{
    \begin{array}{ll}
      \int_{-1}^{1} P_{i-1}(t)\,\D t & 1 \leq i \leq M,\\
      \int_{-1}^{1} P_{i-M-1}(t)\log [(x-t)^{2}+y^{2}]^{1/2}\,\D t
      & M+1 \leq i \leq 2M,\\
      \int_{-1}^{1} P_{i-2M-1}(t)/[(x-t)^{2}+y^{2}]^{1/2}\,\D t & 2M+1
      \leq i \leq 3M,\\
      \int_{-1}^{1} P_{i-3M-1}(t)[(x-t)^{2}+y^{2}]\,\D t & 3M+1 \leq i \leq 4M.
    \end{array}
  \right.
\end{equation}
The system of equations is solved using the appropriate LAPACK
solver~\cite{lapack99}, in the least squares sense when $N>4M$ and in
the minimum norm sense when $N<4M$. The only issue which must be
clarified is the evaluation of the moments, $m_{i}$. In the case when
$y\equiv0$, the integrals are true Cauchy principal values or Hadamard
finite parts and there exist formulae for their evaluation in terms of
associated Legendre functions~\cite{kaya-erdogan87a} or simple finite
part integrals combined with exact Gaussian
quadratures~\cite{carley07a}. In this case, however, such simple
formulae are not available and a different approach is required.

\subsection{Integration of weighted Legendre polynomials}
\label{sec:method:integrals}

The method outlined in \S\ref{sec:rules} for the evaluation of the
quadrature weights $w_{i}$ requires the evaluation of the moments
$m_{i}$ where
\begin{equation}
  \label{equ:moment:def}
  m_{n} = \int_{-1}^{1} u(t) P_{n}(t)\,\D t,
\end{equation}
where weighting function $u(t)$ will be one of
$\log[(x-t)^{2}+y^{2}]^{1/2}$ or $[(x-t)^{2}+y^{2}]^{-n/2}$, $n=1$ or
2. A general method can be applied to finding the integrals of
weighted Legendre polynomials, using basic functional relations and
simple formulae for the integrals of elementary functions, readily
found in standard references~\cite{gradshteyn-ryzhik80}.

Assuming that integrals of the form
\begin{equation}
  J_{n} = \int_{-1}^{1}t^{n} u(t)\,\D t,
\end{equation}
can be evaluated (the formulae required for this paper are given in
the appendix), the expansion of powers of $t$ in terms of Legendre
polynomials~\cite[8.922.1]{gradshteyn-ryzhik80}:
\begin{subequations}
  \label{equ:legendre}
  \begin{eqnarray}
    t^{2n} &= \frac{1}{2n+1}P_{0}(t) + 
    \sum_{k=1}^{n}(4k+1)
    \frac{2n(2n-2)\ldots(2n-2k+2)}{(2n+1)(2n+3)\ldots(2n+2k+1)}
    P_{2k}(t)\\
    t^{2n+1} &= \frac{3}{2n+3}P_{1}(t) + 
    \sum_{k=1}^{n}(4k+3)
    \frac{2n(2n-2)\ldots(2n-2k+2)}{(2n+3)(2n+5)\ldots(2n+2k+3)}
    P_{2k+1}(t)
  \end{eqnarray}
\end{subequations}
can be used to show that:
\begin{subequations}
  \label{equ:legendre:int}
  \begin{eqnarray}
    J_{2n} &=
    \frac{1}{2n+1}\int_{-1}^{1}P_{0}(t)u(t)\,\D t \nonumber\\
    & + \sum_{k=1}^{n}(4k+1)
    \frac{2n(2n-2)\ldots(2n-2k+2)}{(2n+1)(2n+3)\ldots(2n+2k+1)}
    \int_{-1}^{1}P_{2k}(t)u(t)\,\D t,\\
    J_{2n+1} &=
    \frac{3}{2n+1}\int_{-1}^{1}P_{1}(t)u(t)\,\D t \nonumber\\
    & + \sum_{k=1}^{n}(4k+3)
    \frac{2n(2n-2)\ldots(2n-2k+2)}{(2n+3)(2n+5)\ldots(2n+2k+3)}
    \int_{-1}^{1}P_{2k+1}(t)u(t)\,\D t.
  \end{eqnarray}
\end{subequations}
Then, given the integrals $J_{i}$, $i=0,1,\ldots,N$, the corresponding
integrals of the weighted Legendre polynomials, $m_{i}$, can be
evaluated via:
\begin{eqnarray}
  \label{equ:legendre:int:even}
  C_{i}(2n+1)m_{2n} &= 
  J_{2n}-J_{0}/(2n+1)-
  \sum_{j=1}^{n} (4j+1)C_{j}m_{2j},\\
  m_{0} &= J_{0},\quad C_{j} = \frac{2n-2j}{2n+2j+3}C_{j-1},\quad
  C_{0} = \frac{2n}{2n+3}\nonumber
\end{eqnarray}
and
\begin{eqnarray}
  \label{equ:legendre:int:odd}
  D_{i}(2n+3)m_{2n+1} &= 
  J_{2n+1}-3J_{1}/(2n+3)-
  \sum_{j=1}^{n} (4j+3)D_{j}m_{2j+1},\\
  m_{1} &= J_{1},\quad D_{j} = \frac{2n-2j}{2n+2j+5}D_{j-1},\quad
  D_{0} = \frac{2n}{2n+3}.\nonumber
\end{eqnarray}
The evaluation can be carried out in place, with values of $m_{i}$
overwriting $J_{i}$. In double precision, the procedure gives accurate
results for $n\lessapprox32$ before overflow errors cause problems.
For BEM applications, where the shape functions are typically of order
no greater than~3, this causes no special difficulties, but would
limit use of the method in other areas. Obviously, for any application
other than straight elements with $G=\log|\mathbf{x}-\mathbf{x}_{1}|$,
a rule with $M>3$ will be needed, but it is unlikely that $M$ will
need to be greater than~32. 

\subsection{Summary of method}
\label{sec:summary}

To summarize, the procedure for computing a quadrature rule with $N$
points which can integrate functions of the form of
equation~(\ref{equ:quad:func2}) where $a$, $b$, $c$ and $d$ are
polynomials of order up to $M$ is as follows:
\begin{enumerate}
\item find the quadrature points $t_{i}$ for an $N$-point Gaussian
  quadrature, using, for example, the method of Davis and
  Rabinowitz~\cite{davis-rabinowitz75};
\item evaluate the weighted Legendre polynomials $\psi_{ij}$ at each
  of the quadrature points;
\item compute the moments $m_{i}$ using the method of
  \S\ref{sec:method:integrals};
\item solve $[\psi_{ij}]w_{j}=m_{i}$ using an appropriate solver.
\end{enumerate}
The integral of $f(x,y,t)$ is then approximated by:
\begin{equation}
  \int_{-1}^{1}f(x,y,t)\,\D t \approx
  \sum_{i=0}^{N}w_{i}f(x,y,t_{i}). 
\end{equation}

\section{Numerical tests}
\label{sec:tests}

\begin{figure}
  \centering
  \includegraphics{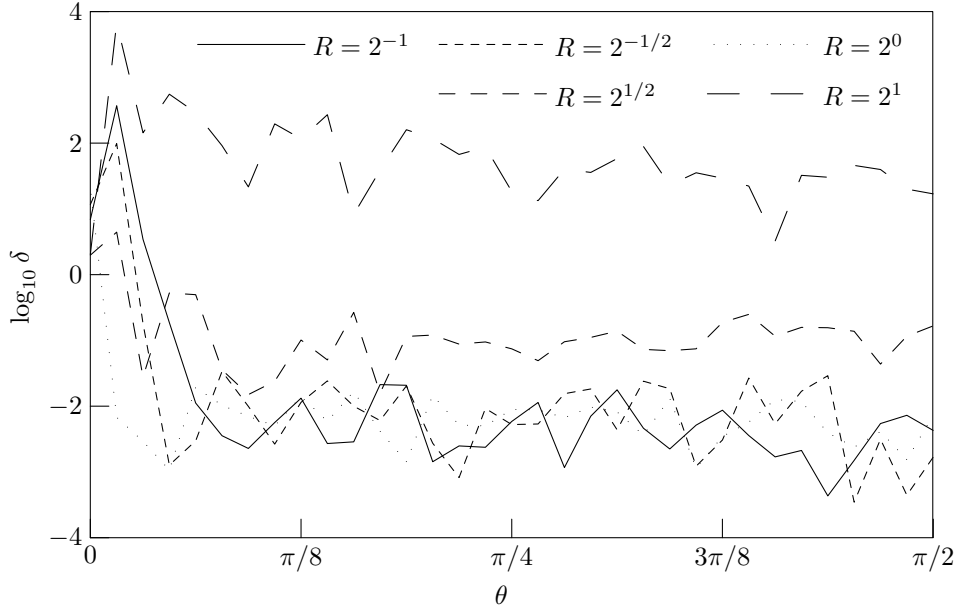}
  \caption{Deviation of quadrature weights from corresponding
    Gauss-Legendre rule: $N=64$, $M=16$.}
  \label{fig:deviation}
\end{figure}

The quadrature method developed in \S\ref{sec:method} is tested by
applying it to the evaluation of a reference integral. Before carrying
out these tests, it is of some interest to examine the behaviour of
the quadrature weights $w_{i}$ with respect to the field point
position. Figure~\ref{fig:deviation} shows the root mean square
difference $\delta$ between the rule of this paper with $N=64$, $M=16$
and a standard~64-point Gaussian quadrature, with
\begin{equation}
  \label{equ:deviation}
  \delta = 
  \left[
    \frac{1}{N}\sum_{i=1}^{N}
    \left(
      w_{i}-v_{i}
    \right)^{2}
  \right]^{1/2},
\end{equation}
where $v_{i}$ are the weights of the standard rule. 

The difference between the rules is shown as a function of $R$ and
$\theta$. It is clear that close to the element ($R$ and/or $\theta$
small), the weights of the new quadrature are very large, as they have
to cancel the large value of the integrand close to the
near-singularity. Further from the element, however, the integrand is
better approximated by a polynomial and the weights come closer to
those of the standard rule. 

\subsection{Accuracy}
\label{sec:tests:accuracy}

The test of accuracy is to examine how well the new quadrature rule
evaluates a reference integral. Table~\ref{tab:error:2} shows the
error in computing:
\begin{equation}
  \label{equ:reference:2}
  I = \int_{-1}^{1}\frac{t^{n}}{(x-t)^{2}+y^{2}}\,\D t,
\end{equation}
for three values of $R$ over the range $\pi/64\leq\theta\leq31\pi/64$.
The case of $\theta=\pi/2$ was ignored because for odd $n$, $I\equiv0$
which would not allow for a meaningful estimate of relative error
which is defined:
\begin{equation}
  \label{equ:epsilon}
  \epsilon = 
  \left[
    \frac{1}{N}
    \sum_{i=0}^{N}
    \frac{(K(\theta_{i})-I(\theta_{i}))^{2}}{I(\theta_{i})^{2}}
  \right]^{1/2},
\end{equation}
where $K$ is the value of $I$ estimated by numerical quadrature. 


Tables~\ref{tab:error:2} and~\ref{tab:error:3} show the error
$\epsilon_{1}$, incurred using the method of this paper, compared to
$\epsilon_{2}$, the error using standard Gauss-Legendre quadrature,
using low and high order rules. The calculation is carried out at
three values of $R$ to examine the change in error with distance from
the element and for various values of $n$. Table~\ref{tab:error:2}
shows the error for $n=0,\dots,3$, the important range for boundary
element calculations, and a rule with $N=16$ and $M=4$.  From the
first two columns of errors, it is clear that the modified rule is far
superior to the standard quadrature: its error is about twelve orders
of magnitude lower than that of the Gaussian quadrature.  Similarly,
for $R=1$, the mean error is much smaller, being no worse than about
$10^{-10}$, rather than $10^{-2}$. Once the field point is far from
the element, however, at $R=2$, both rules have similar accuracy.

\begin{table}
  \centering
  \caption{Root-mean-square error in reference integral computed with
    modified ($\epsilon_{1}$) and standard ($\epsilon_{2}$)
    quadratures at three values of $R$, $N=16$, $M=4$.}
  \label{tab:error:2}
  \begin{tabular}{lllllll}
    \hline
    $R$
    & 
    \multicolumn{2}{c}{$2^{-1}$} &
    \multicolumn{2}{c}{$2^{0}$} &
    \multicolumn{2}{c}{$2^{1}$} \\
    $n$
    &
    \multicolumn{1}{c}{$\epsilon_{1}$} & 
    \multicolumn{1}{c}{$\epsilon_{2}$} & 
    \multicolumn{1}{c}{$\epsilon_{1}$} & 
    \multicolumn{1}{c}{$\epsilon_{2}$} & 
    \multicolumn{1}{c}{$\epsilon_{1}$} & 
    \multicolumn{1}{c}{$\epsilon_{2}$} \\
    \hline
    $0$ & $1.6\times10^{-12}$ & $5.9\times10^{0}$ &
    $3.6\times10^{-11}$ & $1.4\times10^{-2}$ & $2.1\times10^{-16}$ &
    $4.5\times10^{-16}$ \\ 

    $1$ & $2.8\times10^{-13}$ & $3.3\times10^{0}$ & $1.3\times10^{-10}$
    & $1.3\times10^{-2}$ & $1.3\times10^{-16}$ & $7.3\times10^{-16}$ \\ 

    $2$ & $6.0\times10^{-14}$ & $1.9\times10^{0}$ &
    $1.0\times10^{-10}$ & $1.3\times10^{-2}$ & $2.6\times10^{-16}$ &
    $1.5\times10^{-15}$ \\ 

    $3$ & $1.9\times10^{-13}$ & $1.0\times10^{0}$ &
    $9.9\times10^{-11}$ & $1.2\times10^{-2}$ & $6.3\times10^{-16}$ &
    $3.0\times10^{-15}$ \\
   \hline
  \end{tabular}
\end{table}

Table~\ref{tab:error:3} shows similar results for a rule with $N=64$
and $M=16$, compared to data for a standard Gauss-Legendre rule with
$N=64$. Integrals of order up to $n=15$ have been computed and, as
before, at small $R$, the error behaviour of the new rule is orders of
magnitude better than that of the standard technique. At larger $R$,
however, the advantage is not so clear cut: at $R=1$, the error in the
standard rule is around $10^{-12}$, still larger than that from the
method of this paper, but probably acceptable in many applications.
When $R=2$, the Gauss-Legendre rule gives results comparable to those
of the new technique, although it performs better on low order
polynomials. As might be expected, when the distance from the element
is large enough, a high order Gauss-Legendre rule can capture enough
of the behaviour of the integrand to accurately compute the integral.

\begin{table}
  \centering
  \caption{Root-mean-square error in reference integral computed with
    modified ($\epsilon_{1}$) and standard ($\epsilon_{2}$)
    quadratures at three values of $R$, $N=64$, $M=16$.}
  \label{tab:error:3}
  \begin{tabular}{lllllll}
    \hline
    $R$
    & 
    \multicolumn{2}{c}{$2^{-1}$} &
    \multicolumn{2}{c}{$2^{0}$} &
    \multicolumn{2}{c}{$2^{1}$} \\
    $n$
    &
    \multicolumn{1}{c}{$\epsilon_{1}$} & 
    \multicolumn{1}{c}{$\epsilon_{2}$} & 
    \multicolumn{1}{c}{$\epsilon_{1}$} & 
    \multicolumn{1}{c}{$\epsilon_{2}$} & 
    \multicolumn{1}{c}{$\epsilon_{1}$} & 
    \multicolumn{1}{c}{$\epsilon_{2}$} \\
    \hline
$0$ & $1.9\times10^{-10}$ & $1.1\times10^{-01}$ & $6.6\times10^{-15}$ & $5.6\times10^{-12}$ & $1.9\times10^{-12}$ & $3.7\times10^{-16}$ \\
$3$ & $3.4\times10^{-10}$ & $8.4\times10^{-03}$ & $5.4\times10^{-15}$ & $4.2\times10^{-12}$ & $6.5\times10^{-13}$ & $3.0\times10^{-15}$ \\
$6$ & $4.0\times10^{-10}$ & $3.8\times10^{-03}$ & $7.1\times10^{-15}$ & $2.6\times10^{-12}$ & $2.3\times10^{-12}$ & $2.5\times10^{-14}$ \\
$9$ & $4.0\times10^{-10}$ & $8.0\times10^{-04}$ & $8.5\times10^{-15}$ & $9.9\times10^{-13}$ & $1.0\times10^{-12}$ & $2.1\times10^{-13}$ \\
$12$ & $4.0\times10^{-10}$ & $1.4\times10^{-04}$ & $8.4\times10^{-15}$ & $6.6\times10^{-13}$ & $2.8\times10^{-12}$ & $1.7\times10^{-12}$ \\
$15$ & $3.8\times10^{-10}$ & $2.2\times10^{-05}$ & $8.3\times10^{-15}$ & $2.3\times10^{-12}$ & $5.2\times10^{-12}$ & $1.3\times10^{-11}$ \\
\hline



  \end{tabular}
\end{table}

\begin{figure}
  \centering
  \includegraphics{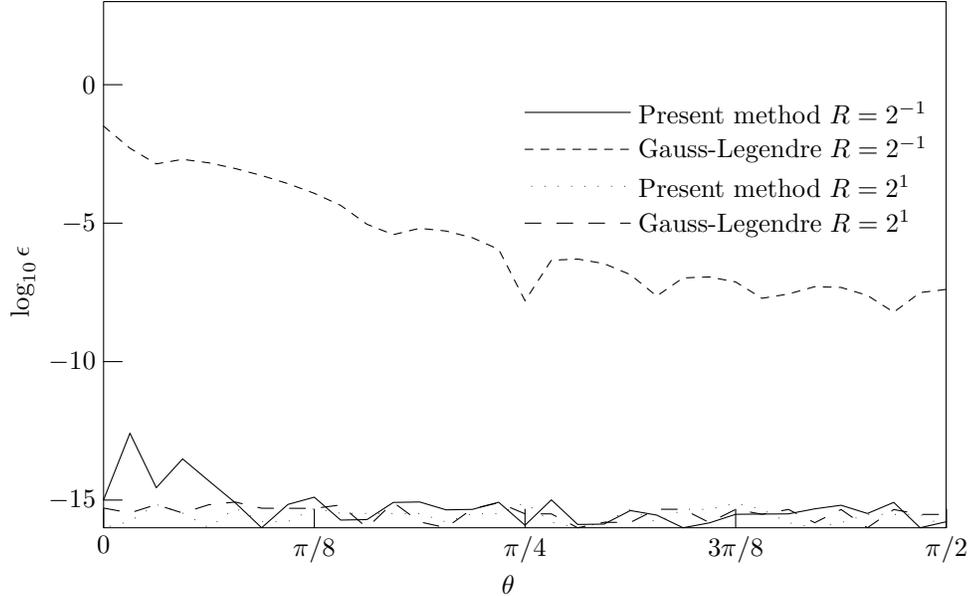}
  \caption{Relative error in logarithmically singular integral, $N=16$,
    $M=4$.}
  \label{fig:error:log:164}
\end{figure}

\begin{figure}
  \centering
  \includegraphics{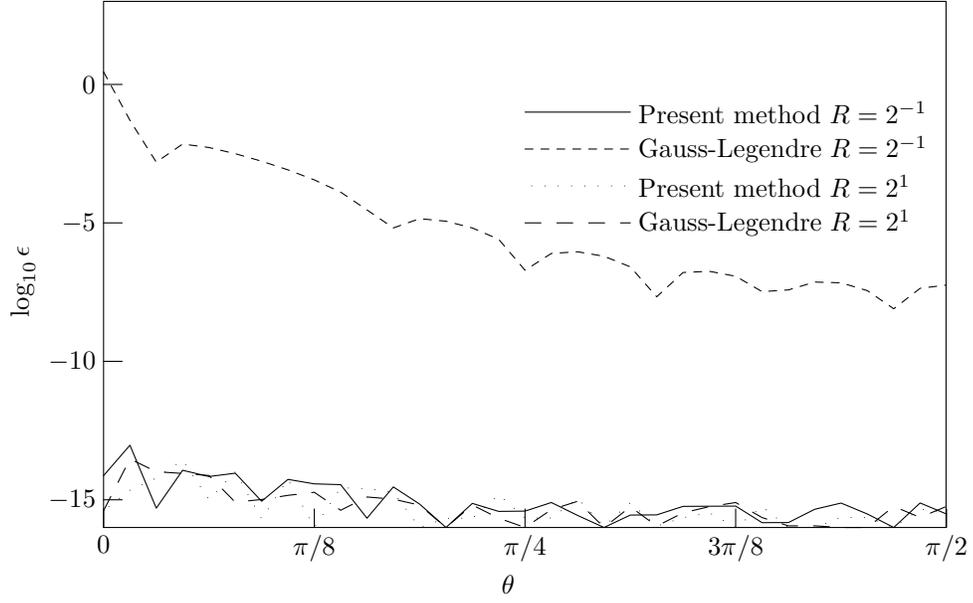}
  \caption{Relative error in near-singular integral, $N=16$, $M=4$.}
  \label{fig:error:1:164}
\end{figure}

\begin{figure}
  \centering
  \includegraphics{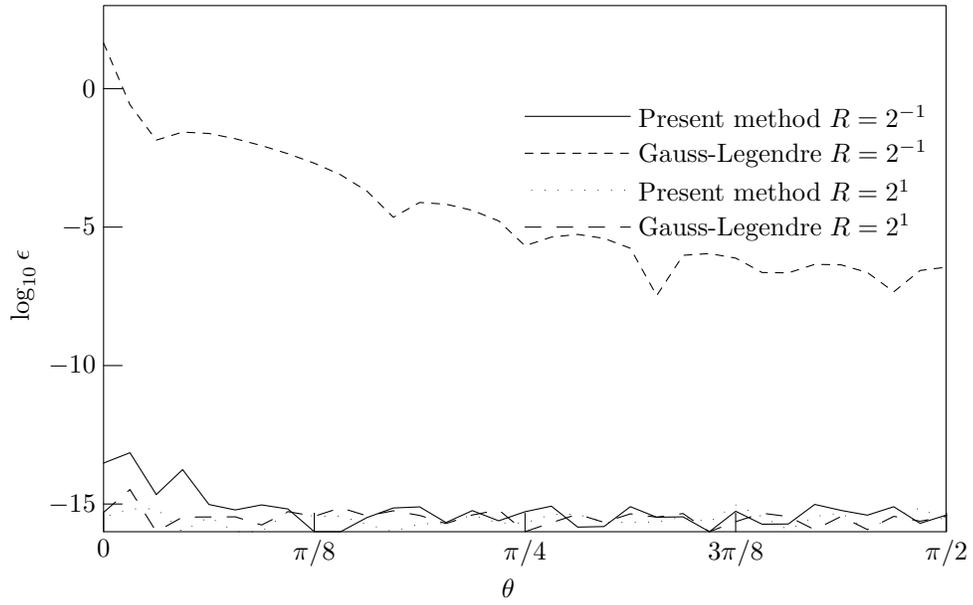}
  \caption{Relative error in near-hypersingular integral, $N=16$, $M=4$.}
  \label{fig:error:2:164}
\end{figure}

To examine the error behaviour in more detail,
figures~\ref{fig:error:log:164}--\ref{fig:error:2:164} show the
relative error in computing:
\begin{align}
  I_{0}^{(\log)}(x,y) &= \int_{-1}^{1} \log[(x-t)^{2}+y^{2}]^{1/2},\\
  I_{0}^{(1)}(x,y) &= \int_{-1}^{1} [(x-t)^{2}+y^{2}]^{-1/2},\\
  I_{0}^{(2)}(x,y) &= \int_{-1}^{1} [(x-t)^{2}+y^{2}]^{-1},
\end{align}
with error $\epsilon$ defined as:
\begin{equation}
  \label{equ:def:error}
  \epsilon = \left|\frac{I_{0}-J_{0}}{I_{0}}\right|,
\end{equation}
where $I_{0}$ is one of $I^{(\log)}_{0}(x,y)$, $I^{(1)}_{0}(x,y)$,
$I^{(2)}_{0}(x,y)$ and $J_{0}$ is the corresponding estimate using the
quadrature rule. In each case, $N=16$ and $M=4$ and a sixteen point
Gaussian quadrature was also applied for comparison. In each case, for
large distances from the element, $R=2$, both quadrature techniques
are accurate, with errors of the order of machine precision. When
$R=2^{-1}$, however, the error incurred using Gaussian quadrature is
unacceptably large, while the modified rule gives very accurate
answers, again of the order of machine precision. As might be
expected, the error from the Gaussian quadrature is smaller as
$\theta\to\pi/2$, due to the greater distance from the element, but it
is never better than about $10^{-6}$, insufficient accuracy for most
applications.

\begin{figure}
  \centering
  \includegraphics{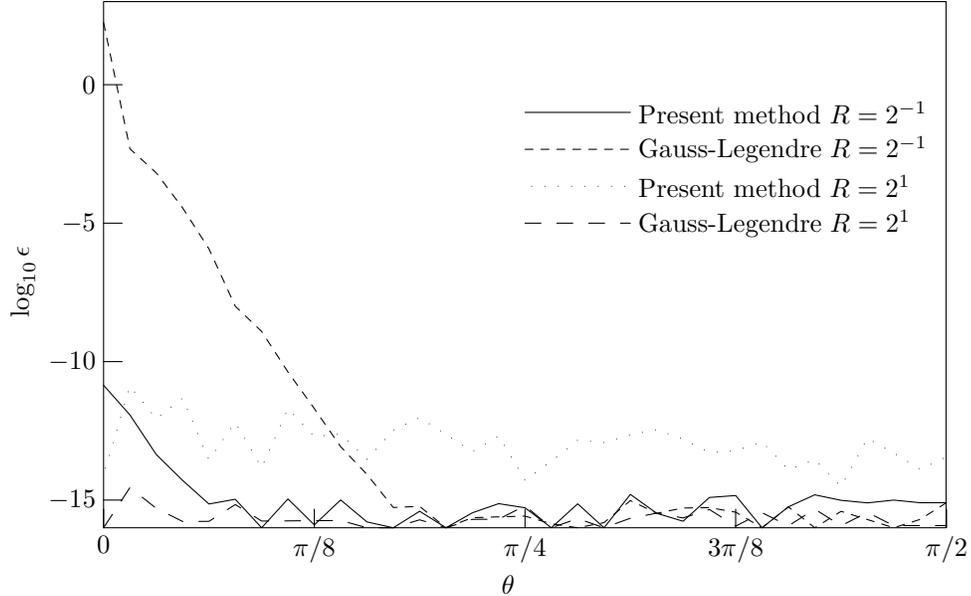}
  \caption{Relative error in near-hypersingular integral, $N=64$,
    $M=16$.}
  \label{fig:error:2:6416}
\end{figure}

Finally, figure~\ref{fig:error:2:6416} illustrates an interesting
point about the error behaviour of the quadrature rule as the number
of quadrature points is increased. It shows the error in the
near-hypersingular integral $I_{0}^{(2)}(x,y)$ evaluated using a rule
with $N=64$ and $M=16$, with the error from a 64-point Gaussian
quadrature shown for comparison. The first point is that the Gaussian
quadrature is able to cope with the singularity for $R=2^{-1}$ when
$\theta\gtrapprox\pi/8$, because it can integrate polynomials of high
enough order to be able to correctly handle the expansion of the
integrand in Legendre polynomials. Also, as in the previous cases, it
can correctly evaluate the integral for $R=2$. 

The modified rule, however, has slightly worse error behaviour than
for $N=4$, with the maximum error being higher for $R=2^{-1}$ and the
error at $R=2$ being greater across the full range of $\theta$. The
error is still small, being less than $10^{-10}$, but the reason for
the increase is unclear. It appears to be due to an ambiguity in
expressing the integrand in terms of Legendre polynomials: the
near-hypersingular part of the integrand $f(t)/[(x-t)^{2}+y^{2}]$ is a
ratio of two polynomials which can be written as a sum of a proper
elementary function and a polynomial. This polynomial term is then
represented twice in the quadrature rule, being handled by the
unweighted Legendre polynomials, $m_{i}$, $1\leq i\leq M$ in
equation~(\ref{equ:moments}) and by the weighted polynomials. An
interesting question for future developments of the technique will be
how best to choose the elementary functions for the quadrature rule to
give optimal accuracy for a given $N$.

\section{Conclusions}
\label{sec:conclusions}

A method of deriving quadrature rules for the evaluation of the
`near-singular' integrals which arise in the boundary element method
has been derived. The performance of the technique has been assessed
by evaluation of reference integrals and it has been found that it
outperforms standard Gaussian quadrature rules with the same number of
nodes for field points close to the element. The error in the integral
increases slightly with the number of points in the rule, a point
which is to be investigated in future work.

\appendix
\section{Integrals of weighted polynomials}
\label{sec:poly:int}

To evaluate the required integrals of Legendre polynomials using the
procedure of \S\ref{sec:method:integrals}, we require a method of
evaluating the integrals:
\begin{subequations}
  \label{equ:base:int}
  \begin{eqnarray}
    I_{n}^{(2)}(x,y) &= 
    \int_{-1}^{1}
    \frac{t^{n}}{(x-t)^{2}+y^{2}}
    \,\D t,\\
    I_{n}^{(1)}(x,y) &= 
    \int_{-1}^{1}
    \frac{t^{n}}{[(x-t)^{2}+y^{2}]^{1/2}}
    \,\D t,\\
    I_{n}^{(\log)}(x,y) &= 
    \int_{-1}^{1}
    t^{n}\log [(x-t)^{2}+y^{2}]^{1/2}
    \,\D t,
  \end{eqnarray}
\end{subequations}
Use of standard formulae~\cite[2.171,2.263,2.728.1]{gradshteyn-ryzhik80}
yields:
\begin{subequations}
  \label{equ:base:form}
  \begin{eqnarray}
    I^{(2)}_{n}(x,y) &= \frac{1 + (-1)^n}{n-1} + 2xI^{(2)}_{n-1}(x,y) - 
    R^{2}I^{(2)}_{n-2}(x,y),\\
    I^{(1)}_{n}(x,y) &= \frac{R_{+}-R_{-}}{n} +
    \frac{2n-1}{n}xI^{(1)}_{n-1}(x,y)
    - \frac{n-1}{n}R^{2}I^{(1)}_{n-2}(x,y),\\
    I^{(\log)}_{n}(x,y) &= \frac{\log R_{+} + (-1)^n\log R_{-}}{2(n+1)}
    - \frac{1}{n+1}I^{(2)}_{n}(x,y) +
    \frac{R\cos\theta}{n+1}I^{(2)}_{n}(x,y),
  \end{eqnarray}
\end{subequations}
where the recursions are seeded with:
\begin{eqnarray*}
  I^{(2)}_{0}(x,y) &= 
  \frac{1}{y}
  \left(
    \tan^{-1}\frac{1-x}{y}
    -
    \tan^{-1}\frac{1+x}{y}
  \right),\quad
  I^{(2)}_{1}(x,y) = \log \frac{R_{+}}{R_{-}} + x  I^{(2)}_{0}(x,y),\\
  I^{(1)}_{0}(x,y) &= \log\frac{R_{+}+1-x}{R_{-}-1-x},
  \quad   
  I^{(1)}_{1}(x,y) = R + x  I^{(1)}_{0}(x,y),
\end{eqnarray*}
and
\begin{eqnarray*}
  R &= (x^{2}+y^{2})^{1/2}, \quad R_{\pm} =
  [(x\mp1)^{2}+y^{2}]^{1/2},\quad \theta = \tan^{-1} y/x.  
\end{eqnarray*}


\end{document}